\documentclass[11pt,a4paper, oneside]{amsart}
\usepackage[utf8]{inputenc}
\usepackage[english]{babel}
\usepackage{amsfonts,amsmath,amssymb,amsthm,enumitem, mathrsfs}
\usepackage[left=2.4 cm,right=2.4 cm,top=2.4 cm,bottom=2.3 cm]{geometry}
\usepackage{setspace}
\usepackage{mathtools}

\usepackage{color, xcolor}
\definecolor{LinkBlue}{RGB}{0,92,170}
\definecolor{CiteGreen}{RGB}{0,128,90}
\definecolor{UrlPurple}{RGB}{120,0,160}
\PassOptionsToPackage{
  unicode,
  colorlinks=true,
  linkcolor=LinkBlue,
  citecolor=CiteGreen,
  urlcolor=UrlPurple,
  pdfborder={0 0 0}
}{hyperref}
\usepackage{ifpdf}
\ifpdf
    \usepackage[pdftex]{graphicx}
    \DeclareGraphicsExtensions{.png,.pdf,.jpg}
\else
   \usepackage[dvips]{graphicx}
   \DeclareGraphicsExtensions{.eps}
\fi
\graphicspath{{.}{figures/}}
\numberwithin{equation}{section}

\usepackage{doi}
\hypersetup{
  colorlinks=true,
  linkcolor=LinkBlue,
  citecolor=CiteGreen,
  urlcolor=UrlPurple,
  pdfborder={0 0 0}
}
\usepackage{comment}

\newtheorem{theorem}{Theorem}[section]
\newtheorem{lemma}[theorem]{Lemma}

\newtheorem*{theorem*}{Theorem}
\newtheorem*{lemma*}{Lemma}
\newtheorem*{proposition*}{Proposition}
\newtheorem*{corollary*}{Corollary}

\renewcommand\tilde{\widetilde}

\newcommand{\dkl}{\mathrm{d}_{\mathrm{KL}}}

\renewcommand{\phi}{\varphi}

\def\1{\mathbf{1}}

\def\XXint#1#2#3{{\setbox0=\hbox{$#1{#2#3}{\int}$ }
\vcenter{\hbox{$#2#3$ }}\kern-.57\wd0}}

\def\W{\mathscr{W}}
\begin{document}

\title{The Wasserstein cost of the Girsanov theorem}
\author{Simon Coste}
\address{Laboratoire de Probabilités, Statistique et Modélisation, Universit\'e Paris Cit\'e, Paris, France}
\email{simon.coste@u-paris.fr}
\date{\today}

\begin{abstract}
This is an expository note on the asymptotics of the Wasserstein distance between the empirical measure of infinite-dimensional random paths and their common distribution. We present the classical proof for these asymptotics in the case of the Wiener measure, unifying quantization arguments from Dereich and coauthors \cite{dereich2008coding,dereich2005randomly} and Boissard-Le Gouic \cite{boissard2014mean}, and we adapt it for Itô diffusions sampled through the Girsanov theorem.

More precisely, let $\mathbb{Q}$ be the law of the solution of a stochastic differential equation with a bounded Lipschitz drift. Let $W_1,\ldots,W_n$ be independent Brownian paths, and let $Z_i$ be their Girsanov weights for $\mathbb{Q}$. The reweighted empirical measure
$$\widehat\mu_n=\frac{\sum_{i=1}^n Z_i\delta_{W_i}}{\sum_{i=1}^n Z_i}$$
approximates $\mathbb{Q}$ when $n$ is large. We prove that the Wasserstein cost $\W_p^p(\widehat\mu_n,\mathbb{Q})$ is of order $(\log n)^{-p/2}$, in expectation and almost surely, for every fixed $p\geqslant1$. 

In the Brownian case where $Z_i=1$ these results show that $\sqrt{\log n}\W_p(\widehat\mu_n,\mathbb{Q})$ is asymptotically sandwiched between $\sqrt{2\lambda_d}$ and $4\sqrt{2\lambda_d}$, where $\lambda_d$ is the first eigenvalue of the Laplacian operator in the unit ball. We provide numerical experiments to estimate the (conjectured) limiting constant in dimension $d=1$. 
\end{abstract}

\maketitle
\setcounter{tocdepth}{1}
\tableofcontents

\section{Introduction and statement of the results}

\subsection{Importance sampling} \label{sec:is}

Let $f,g$ be two probability densities on $\mathbb{R}^d$, with the same support. For any bounded measurable function $\varphi$, the identity
\begin{equation}\label{eq:is_identity}
    \mathbb{E}_{X\sim g}[\varphi(X)]=\int \varphi g=\int \varphi\frac{g}{f}f = \mathbb{E}_{Y\sim f}\left[\varphi(Y)\frac{g(Y)}{f(Y)}\right]
\end{equation}
suggests a simple way to approximate expectations of functions of a random variable with respect to $g$: we first draw samples $Y_1,\ldots,Y_n$ from the density $f$ (called the \emph{proposal}), and then we reweight each sample $Y_i$ by $w(Y_i)$, where $w(y) = g(y)/f(y)$ is the \emph{importance ratio}. Due to the law of large numbers, the empirical mean $n^{-1}\sum_{i=1}^n \varphi(Y_i)w(Y_i)$ converges to $\mathbb{E}[\varphi(Y)w(Y)] = \mathbb{E}_{X\sim g}[\varphi(X)]$ when $n$ is large. To ensure that the estimator for the constant function 1 is always equal to 1, it is even better to use the \emph{self-normalized} version of this estimator, which is 
\begin{equation}\label{eq:snis}\frac{\sum_{i=1}^n w(Y_i)\varphi(Y_i)}{\sum_{i=1}^n w(Y_i)}.
\end{equation}
This technique is known as \emph{importance sampling}, and works well in contexts where $f$ is easy to sample from. Using \eqref{eq:snis} also has a crucial advantage: it works when the densities $f,g$ are only known up to constants (for example, in Monte Carlo integration and Bayesian computation; see \cite{kahn1953methods,agapiou2017importance,elvira2021advances}). 
However, it is known to practitioners that this technique comes with a caveat: sometimes, the variance of the weights can be too large, leading to poor estimators with a slow decay in the number of samples $n$. So, given a predefined acceptable error, how many samples are needed to reach this error? 

The landmark paper by Chatterjee and Diaconis \cite{chatterjee_diaconis} describes an ``information-theoretic threshold'' at $n_*= \exp(\dkl(g\mid f))$: taking more samples than $n_*$ is sufficient to reach a low error \emph{for a fixed observable $\varphi$}, but fewer than $n_*$ can lead to very poor estimators \emph{for some pathological $\varphi$}, see also \cite{coste2023importance} for an accessible exposition. In our previous work \cite{coste2026wasserstein}, we instead considered the approximation of all the $\varphi$ simultaneously, i.e. the whole measure $$\hat{g}_n = \frac{\sum_{i=1}^n w(Y_i)\delta_{Y_i}}{\sum_{i=1}^n w(Y_i)}, $$ and not just the expectation of a fixed observable. The right metric here is the Wasserstein distance. Under some assumptions and in dimension $d\geqslant3$, we proved that there are two constants $0<C_1 \leqslant C_2 < \infty$ depending on $p$ and $d$ such that for $n$ large enough,
\begin{equation}\label{eq:previous}
    C_1 \int gf^{-p/d}\leqslant \liminf \& \limsup ~ n^{p/d}\mathbb{E}\W_p^p(\hat g_n,g) \leqslant C_2    \int gf^{-p/d}.   
\end{equation}
When $p=2$ we know that $C_1 = C_2$, hence this is a true limit. The rate of decay $n^{-p/d}$ is the classical one for empirical measures, and we refer to \cite{AKT84,Ta94,BaBo13,FoGu15,goldman2024optimal} and the discussion in \cite{coste2026wasserstein} for further information on this important problem. The main point of our work was the extra term $\int gf^{-p/d}$, explaining how the asymptotics depend on the interaction between the target $g$ and the proposal $f$. The exponent $1/d$ is often interpreted as an instance of the curse of dimensionality phenomenon: when $d$ is high, the convergence is extraordinarily slow. In this note, we revisit these results, but on infinite-dimensional spaces. A popular instance of the importance sampling technique is known under the name of the \emph{Girsanov theorem}, which gives the density of the law of certain Itô processes with respect to the Wiener measure. This is a classical tool in stochastic calculus and its applications, notably in mathematical finance, see \cite{newton1994variance,glasserman1999asymptotically,glasserman2004monte}. 

\bigskip 

\noindent\textbf{Contribution.} The goal of this short note is to derive results similar to \cite{coste2026wasserstein} for the Girsanov theorem. The proof is an easy adaptation of the Wasserstein asymptotics for the Wiener measure. These asymptotics are implicitly known, but the proof was in various papers, notably those by Dereich and coautors (see \cite{dereich2008coding,dereich2006high,dereich2003link,dereich2005randomly}) and Boissard-Le Gouic \cite{boissard2014mean}. We hope that this note will give a self-contained, easy exposition of the Wiener case, in addition to the new result on the Girsanov theorem, and additional numerics on the unknown limiting constants therein.

\subsection{Notations}In the rest of this note, the parameter $p\geqslant 1$ and the dimension $d\geqslant 1$ are fixed. For two nonzero sequences $(a_n)$ and $(b_n)$, we write $a_n = O(b_n)$ if there is a constant $C>0$ such that $|a_n| \leqslant C|b_n|$ for all $n$ sufficiently large. The constant $C$ may depend on $p$ and $d$, but not on $n$. Similarly we write $a_n = o(b_n)$ if the ratio $|a_n|/|b_n|$ tends to 0 as $n$ tends to infinity and $a_n \sim b_n$ if the ratio $a_n / b_n$ tends to 1.

\subsection{Brownian paths and Wasserstein distance}\label{sec:brownian}

We denote by $W$ and $W_i$ ($i\in \mathbb{N}^*$) independent standard Brownian motions in $\mathbb{R}^d$. We reserve the letter $X$ for the target Itô process. These processes are seen as random variables on the space
$$\mathscr{C}=\{ f : [0,1] \to \mathbb{R}^d \text{ continuous, with } f(0)=0 \},$$
endowed with the uniform topology induced by the norm $\Vert f \Vert = \sup_{t \in [0,1]}|f(t)|$, where $|\cdot|$ is the usual Euclidean norm. We denote by $\mathbb{W}$ the Wiener measure on this space. The empirical distribution of the Brownian paths is
\begin{equation}\label{eq:brownian_empirical}
    \mu_n=\frac1n\sum_{i=1}^n\delta_{W_i}.
\end{equation}
Given two probability measures $\mu,\nu$ with finite $p$th moments, their Wasserstein distance is defined by
\begin{equation}\label{eq:wasserstein}
    \W_p^p(\mu,\nu)=\inf_{\pi\in\Gamma(\mu,\nu)}\int\Vert f-g\Vert^p\,\pi(df,dg),
\end{equation}
where $\Gamma(\mu,\nu)$ is the set of couplings of $\mu$ and $\nu$, i.e. pairs of random variables $(F,G)$ such that $F$ has law $\mu$ and $G$ has law $\nu$, but with no constraints whatsoever on their joint distribution. 

Let $\lambda_d$ be the first eigenvalue of the operator $-\tfrac12\Delta$ in the Euclidean unit ball; here $\Delta$ is the Laplacian operator. In particular, $\lambda_1=\pi^2/8$ as we recall in Appendix \ref{appendix:small_balls}.

\begin{theorem}\label{thm:brownian}
For every $p\geqslant1$, almost surely one has
\begin{equation}\label{eq:brownian_rate}
\begin{aligned} 
\liminf_{n\to\infty}\sqrt{\log n}\W_p(\mu_n,\mathbb{W})&\geqslant\sqrt{2\lambda_d},\\
\limsup_{n\to\infty}\sqrt{\log n} \W_p(\mu_n,\mathbb{W})&\leqslant4\sqrt{2\lambda_d}.
\end{aligned}
\end{equation}
\end{theorem}

\subsection{Girsanov reweighting}\label{sec:girsanov}

Let $h:[0,1]\times\mathbb{R}^d\to\mathbb{R}^d$ be bounded, continuous, globally Lipschitz in space, and uniformly Lipschitz in time. Let $X$ be the unique solution of the stochastic differential equation (SDE)
\begin{equation}\label{sde}
    dX_t=h(t,X_t)dt+dW(t),\qquad X_0=0,
\end{equation}
and denote by $\mathbb{Q}$ its distribution on $\mathscr{C}$. The Girsanov theorem says that $\mathbb{Q}$ is absolutely continuous with respect to $\mathbb{W}$, and that its Radon--Nikodym derivative with respect to a Brownian sample $W_i$ is given by
\begin{equation}\label{eq:girsanov_weight}
    Z_i=\exp\left(\int_0^1 h(t,W_i(t)) dW_i(t) -\frac12\int_0^1|h(t,W_i(t))|^2dt\right).
\end{equation}
In other words, for every bounded measurable functional $\varphi$, we have
\begin{equation}\label{eq:girsanov_expectation}\mathbb{E}[\varphi(X)]=\mathbb{E}[Z_i\varphi(W_i)].\end{equation}
Consequently, we approximate $\mathbb{Q}$ by the weighted empirical measure
\begin{equation}\label{eq:weighted_empirical}
\widehat\mu_n=\frac{\sum_{i=1}^n Z_i\delta_{W_i}}{nS_n},
\end{equation}
where $S_n=n^{-1}\sum_{i=1}^n Z_i$. The denominator ensures that $\widehat\mu_n$ is a probability measure. The following theorem shows that the Brownian bounds of the preceding theorem still hold.

\begin{theorem} \label{thm:girsanov}
For every $p\geqslant1$, almost surely one has 
\begin{equation}\label{eq:girsanov_rate}
\begin{aligned}  
\liminf_{n\to\infty}\sqrt{\log n}\,\W_p(\widehat\mu_n,\mathbb{Q})&\geqslant\sqrt{2\lambda_d},\\
\limsup_{n\to\infty}\sqrt{\log n}\,\W_p(\widehat\mu_n,\mathbb{Q})&\leqslant4\sqrt{2\lambda_d}.
\end{aligned}
\end{equation}
\end{theorem}

Of course, these bounds do not match. It is thus possible that the true asymptotics might be of the form $\W_p(\hat{\mu}_n, \mathbb{Q})\sim c /\sqrt{\log n}$, where $c$ depends on $\mathbb{Q}$ in a non-trivial way, as in the paper \cite{coste2026wasserstein}. We were not able to go in this direction.

\begin{figure}[tbp]
    \centering
    \includegraphics[width=.93\linewidth]{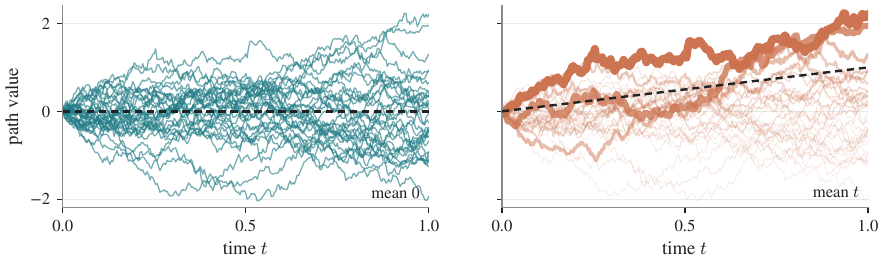}
    
    \caption{The same $n=32$ realizations of Brownian paths before and after reweighting by $Z_i=e^{W_i(1)-1/2}$. Line widths are proportional to the normalized masses, and on the right so is the opacity. The dashed lines are the population means, chosen to be $0$ and $t$.}
    \label{fig:weights}
    \end{figure}

\subsection{Numerical experiments in one dimension }\label{sec:numerics}

In dimension $d=1$, Theorems \ref{thm:brownian}-\ref{thm:girsanov} leave a big gap between the two bounds, $\pi/2\approx 1.57$ and $2\pi\approx 6.28$. Supposing that the quantity $\sqrt{\log n}\W_p(\mu_n, \mathbb{W})$ is converging (which is purely conjectural, but nevertheless reasonnable), where should the limiting constant lie? We explore this question numerically\footnote{The code and saved results are available in \cite{ghrepo}. } for Brownian samples. For simplicity we only investigate the case $p=1$, $d=1$. Directly computing $\W_1(\mu_n, \mathbb{W})$ is impossible, so we have to resort to non-rigorous tricks and proxies. 

\subsubsection{Using the proxy $L_n$}Any coupling of $\mu_n$ with $\mathbb W$ must send each path at least as far as its nearest sample, so we always have
$$L_n=\mathbb E_W\left[\min_{i\leqslant n} \Vert W-W_i\Vert\right] \leqslant\W_1(\mu_n,\mathbb W). $$
Here $\mathbb E_W$ only averages over the new path $W$. Dereich and Lifshits \cite[Theorems~5.1 and~6.1]{dereich2005randomly} show that $\sqrt{\log n}L_n$ converges toward a constant $\sqrt{\kappa}$: they show that if $W$ is an independent Brownian path, then
\begin{equation}\label{eq:dereich}-\varepsilon^2\log\mathbb W(B(W,\varepsilon))\longrightarrow \kappa\quad\text{almost surely, and }
 \mathbb E[L_n]\sim\sqrt{\kappa/\log n}.\end{equation}
Our proof for Theorem \ref{thm:brownian} actually reproduces a part of their work, namely that $\liminf \sqrt{\log n}L_n$ is greater than $\sqrt{2\lambda_d}$, which shows that $\kappa$ is itself greater or equal to $2\lambda_d$.

However, sending $W$ toward its closest sample does not give a coupling between $\mu_n$ and $\mathbb W$. The mass assigned by this coupling to the sample $W_i$ is $p^{(n)}_i = \mathbb{W}(C_i)$, the Wiener mass of the Voronoï cell of $W_i$ among the $n$ samples $(W_i)$. In general, $p^{(n)}_i $ is not equal to $1/n$. We define 
$$\eta_n = \sum_{i=1}^n p^{(n)}_i \delta_{W_i}.$$
Clearly, $L_n$ is exactly equal to the cost $\W_1(\eta_n,\mathbb{W})$, so by the triangle inequality
\begin{equation}\label{eq:boundLn}L_n\leqslant \W_1(\mu_n,\mathbb W) \leqslant L_n+\W_1(\mu_n,\mathbb \eta_n).\end{equation}

\begin{figure}[tbp]
    \centering
    \includegraphics[width=\linewidth]{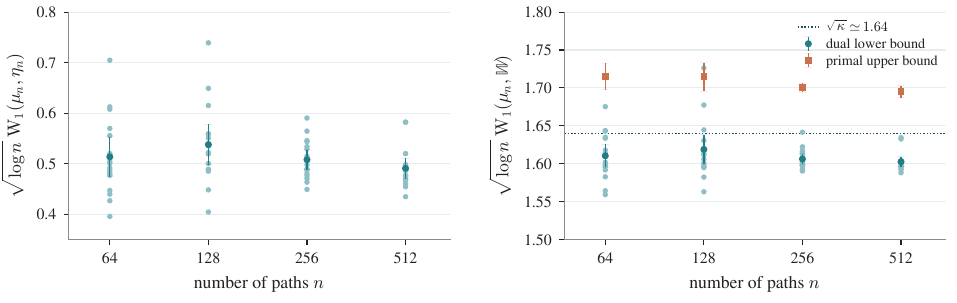}
    \caption{Left: normalized cost of transporting $\mu_n$ onto $\eta_n$, for 16 independent samples at each $n$. Right: bounds on the full normalized cost $\sqrt{\log n}\,\W_1(\mu_n,\mathbb W)$ for the same samples, on the $2^{14}$ grid. The blue points are the dual lower bounds and the brown squares are the primal upper bounds; the dotted line is $\sqrt\kappa\simeq 1.64$. The bars are twice the standard error of the mean. }
    \label{fig:voronoi-transport}
\end{figure}

To estimate $\kappa$, we can either use the direct definition, or use \eqref{eq:dereich} and compute the random small-ball probability limit on the LHS. Computing these probabilities for $24$ independent paths and extrapolating as the time step decreases gives\footnote{It can be noted that the value $\Gamma(1/3)\approx 2.6789385...$ is notably close to 2.68, but we see no reason why the real constant $\kappa$ should coincide with $\Gamma(1/3)$. } $\kappa\simeq2.68$, hence $\sqrt \kappa\simeq1.64$. The statistical error on $\kappa$ is about $0.02$, and the extrapolation adds uncertainty of order $0.1$. But then, is the second term in \eqref{eq:boundLn} negligible? We were suspecting that \emph{yes}, when in reality our numerics could not allow to conclude. To check if $\W_1(\mu_n,\mathbb \eta_n)$ is negligible compared to $1/\sqrt{\log n}$, we computed it for $n=64,128,256,512$, with the Voronoï masses $p^{(n)}_i$ estimated from $100n$ fresh paths, and we found
$$\sqrt{\log n}\,\W_1(\mu_n,\eta_n)\simeq 0.51,\; 0.54,\; 0.51,\; 0.49, $$
which does not seem at all to go to 0, see Figure \ref{fig:voronoi-transport}.

\begin{figure}[tbp]
    \centering
    \includegraphics[width=\linewidth]{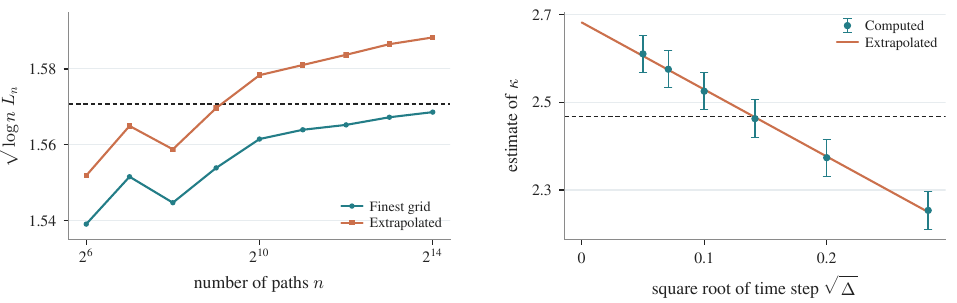}
    \caption{Estimation of the nearest-sample constant. Left: $\sqrt{\log n}L_n$, averaged over 8 samples, with $200$ fresh paths per sample. The blue curve uses a grid with $2^{14}$ points; the brown curve extrapolates from three time grids. Right: estimates of $\kappa$ as the time step $\Delta$ decreases. The bars show twice the standard error over $24$ paths; the line is a fit in $\sqrt\Delta$, with intercept $2.68$. 
    The dashed levels are $\pi/2$ (left) and $\pi^2/4$ (right). }
    \label{fig:brownian-numerics}
\end{figure}

\subsubsection{Bounds on the Wasserstein distance}Since we cannot easily confirm that $\W_1(\mu_n,\mathbb W)$ has the same asymptotics as $L_n$, we tried to approximate directly the latter. To do so, we use the Kantorovich duality. For every $v\in\mathbb R^n$ (the dual potentials), the dual objective function is
$$F(v)=\frac1n\sum_{i=1}^n v_i+\mathbb E_W[\min_{i\leqslant n}(\Vert W-W_i\Vert-v_i)]\leqslant \W_1(\mu_n,\mathbb W).$$
Kantorovich duality ensures that $\sup_v F(v) = \W_1(\mu_n,\mathbb{W})$. We did this for the same 16 samples at $n=64,128,256,512$ as above. The expectation inside $F$ is fitted on $100n$ fresh paths, and then we optimize in $v$.  

The results are in Figure~\ref{fig:voronoi-transport} (right). On the $2^{14}$ grid, the mean of $\sqrt{\log n}\,F(v)$ is $1.610$, $1.619$, $1.606$ and $1.602$ at $n=64,128,256,512$, with a standard error of about $0.01$ at the two smaller sizes and $0.003$ at $n=512$. The corresponding upper bounds are $1.715$, $1.715$, $1.700$ and $1.695$. In particular, the value $\sqrt\kappa\simeq 1.64$ is always inside the confidence interval. 

The comparison with the proxy $L_n$ might be more interesting: indeed, $\sqrt{\log n}\,L_n$ goes from $1.539$ to $1.558$, so $\sqrt{\log n}\,(F(v)-L_n)$ is $0.071$, $0.072$, $0.053$ and $0.044$. At $n=512$ the excess of $\W_1$ over $L_n$ is between $0.044 / \sqrt{\log n}$ and $0.137 /\sqrt{\log n}$, which is less than the $0.49/\sqrt{\log n}$ needed to rebalance the weights of the atoms of $\eta_n$. Most of the cost of $\W_1(\mu_n,\eta_n)$ is therefore an artefact of forcing the nearest-sample coupling to have uniform masses: this suggests that there must be a much better coupling than the "closest point" one. 

\subsubsection{Conclusion}This is globally disappointing, because we could not really make predictions on the (conjectured) limit of $\sqrt{\log n}\W_1(\mu_n, \mathbb{W})$ and it is not clearer that this limit is different than $\sqrt{\kappa}$, as the human author of this note originally suspected. At the accessible sizes the full Wasserstein cost seems to be slightly above $L_n$. The candidate constant $\sqrt{\kappa}\simeq 1.64$ remains plausible, our computations do not confirm it nor exclude it.

\section{Discussion}

\subsection{Lean formalization}We also provide a Lean 4 formalization of the results, available at the public repository \cite{ghrepo}. Since \texttt{Mathlib} lacks tools in stochastic calculus, we rely on the \texttt{MathFin} library for the Girsanov theorem. The formalization was done in around 6 hours of wall-clock time by GPT-Astra-6. The concordance between the Lean formalized results and the results written in this paper was audited by Fable 5.1 and GPT-Astra-6. 

\subsection{Limitations}The Girsanov theorem is often used to compute expectations of observables with respect to Itô processes, especially in mathematical finance. However, the reason for using Girsanov’s theorem with Brownian samples and not directly sampling the Itô process is often to reduce the variance of the estimators. 

Typically, suppose for simplicity that the target process is a Brownian motion $W$, and one wants to integrate an observable $\varphi$ whose support is on extreme trajectories of $W$. In the vanilla Monte-Carlo estimate $n^{-1}\sum_{i=1}^n\varphi(W_i)$, most terms would be zero, although the estimator is in principle unbiased. For this reason, it is easier to sample Itô processes $X_i$ that are designed to take values in the support of $\varphi$ (for example with a drift), and then to tilt them back toward the Wiener measure with reweighting, to compensate for the different distribution, leading to estimators of the form $\sum Z_i \varphi(X_i) / \sum Z_i$. Since most of the paths are in the support of $\varphi$, this can effectively reduce the variance of the estimator. In this view, our result is not very informative, because it gives a worst-case bound, and not a bound on a specific observable of interest (say, $\varphi$). Such bounds could be found with more classical tools. 

Moreover, we only studied the Girsanov theorem with Brownian proposals and Itô targets. What happens when the proposal is not itself Brownian is not clear.

\subsection{Related work}\label{sec:related_work}

The Brownian rate $(\log n)^{-1/2}$ is classical. Dereich et al. \cite{dereich2003link} relate Gaussian small-ball probabilities to deterministic and random quantization; our use of auxiliary Brownian paths to obtain a deterministic codebook belongs to this  approach. Sharp quantization asymptotics are available for Brownian motion \cite{dereich2006high} and diffusions \cite{dereich2008coding}. These concern optimally chosen supports and masses, whereas our empirical approximation must also control the randomness in the points. Our lower bound uses the actual random support: $W-W_i\stackrel{d}{=}\sqrt2W$ gives the constant $\sqrt{2\lambda_d}$.

Boissard and Le Gouic \cite[Theorems~1.1 and~1.4]{boissard2014mean} bound empirical Wasserstein distances through covering numbers and cell-mass fluctuations; Gaussian isoperimetry then gives the Brownian rate in $\W_2$, with concentration. Weed and Bach \cite[Proposition~1 and Section~4]{weed2019sharp} refine this multiscale method using coverings of sets of large measure and Wasserstein dimensions. Their theory only applies to compact spaces, but cannot apply to our setting on the (unbounded) Wiener space. 

The closest finite-codebook argument is He and Yuan \cite[Proposition~4.6]{he2026sharp}, which was put online during the writing this draft. They prove a result similar to Theorem \ref{thm:brownian} (that is, no Girsanov reweighting) but for general Itô processes, even with random, path-dependent coefficients (but under moment assumptions). Our Brownian upper proof has a similar structure as theirs, with Gaussian small-ball estimates keeping constants explicit, but our techniques give sharper constants. 

We are not aware of any work on measuring the accuracy of the Girsanov theorem (or any infinite-dimensional reweighting scheme) in Wasserstein distance.

\subsection{AI disclosure}The research program of this paper was initiated by the human author. GPT-Astra-6 extracted several arguments from the litterature and gave the elements needed to adapt them to the Girsanov case. The writing of this paper was entirely delegated to the human author, except for Subsection \ref{sec:related_work} and the appendices (GPT-Astra-6). The code for the numerical simulations was entirely written by Grok-4.7. 

\subsection{Acknowledgments}This work was carried out while its human author was visiting the Centre de Mathématiques Appliquées at École Polytechnique. The human author thanks another human, Michael Goldman, for useful discussions on the topic.

\section{A simple exposition of the Wiener quantization argument}

In this section, we explain the crucial construction of a \emph{quantization} of the Wiener measure, that is, a finite family of functions such that most of the Wiener mass is very close to one of these functions, with precise quantitative controls. There were many works on Wiener quantization, see for example \cite{dereich2005randomly,graf2000foundations,dereich2008coding}. We extract a simple proof designed for our needs. Most arguments here are classical. We can summarize the main result of this section as follows. 

\begin{lemma}[Wiener Quantization]
Given any $\varepsilon>0$ and $\gamma>2\lambda_d$, there is a finite family of $m\leqslant e^{\gamma \varepsilon^{-2}}$ functions $f_1, \dotsc, f_m \in \mathscr{C}$ with norm smaller than $2\lambda_d \varepsilon^{-1}(1+o(1))$, and such that with probability at least $1 - 2e^{-O(\varepsilon^{-2})}$, a Brownian path $W$ is at distance at most $2\varepsilon$ from one of these functions.
\end{lemma}

\subsection{Small-ball probabilities}\label{sec:small_balls}

For any element $f\in\mathscr{C}$ and radius $r>0$, we denote by
$$B(f,r)=\{g\in\mathscr{C}:\Vert g-f\Vert\leqslant r\}$$ the closed ball around $f$ in the uniform norm. Its volume under the Wiener measure,
\begin{equation}\label{eq:small_ball_definition}
    b_d(r)=\mathbb{W}(B(0,r))=\mathbb{P}(\Vert W\Vert\leqslant r),
\end{equation}
is called the Brownian small-ball probability: it represents the probability that a Brownian path stays in a tubular neighborhood of radius $r$ around the origin. The small-ball probability has been thoroughly studied in the literature. We will use two classical facts: first, when $r \to 0$, 
\begin{equation}\label{eq:bound_bd}
    -\log b_d(r)=\frac{\lambda_d}{r^2}+O(1).
\end{equation}
Second, Anderson's inequality states that for any fixed $f$,
\begin{equation}\label{eq:anderson}
    \mathbb{W}(B(f,r))\leqslant b_d(r).
\end{equation}
In other words, the ball of radius $r$ in $\mathscr{C}$ with the largest Wiener volume is the ball centered at zero. In Appendix~\ref{appendix:small_balls}, we give extra details on these facts, like the explicit formula for $b_1$.

\subsection{General bound on the small-ball probabilities}\label{sec:epsilon_net}

We start by proving a general lower bound on the balls $B(f,r)$, but for functions $f$ that are not necessarily zero. Let $H$ be the \emph{Cameron--Martin space}, i.e. the space of absolutely continuous paths $f$ starting at zero, with the so-called $H^1$-norm
$$\Vert f\Vert_H^2=\int_0^1|f'(t)|^2\,dt<\infty.$$
We denote by $H_R=\{f\in H:\Vert f\Vert_H\leqslant R\}$ the ball of radius $R$ induced by this norm. It should be noted that the Cauchy--Schwarz inequality implies that the uniform norm is smaller than the $H^1$ norm, $\Vert f\Vert\leqslant\Vert f\Vert_H$. 

Let $f$ be in the Cameron--Martin space and consider the process $W(t) + f(t)$. It is obviously a solution of the SDE $dX_t = f'(t)dt + dW(t)$, hence Girsanov's theorem applies, and gives that for any functional $\Psi:\mathscr{C}\to\mathbb{R}_+$, 
$$\mathbb{E}[\Psi(W+f)]=\mathbb{E}\left[e^{\int_0^1 f'(t)dW_t - \frac{1}{2} \int_0^1 |f'(t)|^2 dt} \Psi(W)\right].$$
This specific instance of the Girsanov theorem is often known under the name of \emph{Cameron--Martin formula}. With $\Psi(x) = \mathbf{1}_{\Vert x\Vert\leqslant u}$, we get that
$$\mathbb{W}(B(f,u))=e^{-\Vert f\Vert_H^2/2}\mathbb{E}\left[e^{-\int_0^1 f'(t)dW(t) }\mathbf{1}_{W \in B(0,u)}\right].$$
Since the distribution of $W$ is symmetric, one can replace $W$ by $-W$ in the expectation, and take the average, leading to 

\begin{align*}
\mathbb{W}(B(f,u))
&=e^{-\Vert f\Vert_H^2/2}\mathbb{E}\left[\mathbf{1}_{\Vert W\Vert\leqslant u}\cosh\left(\int_0^1 f'(t) dW(t)\right)\right]\\
&\geqslant e^{-\Vert f\Vert_H^2/2}b_d(u).
\end{align*}
The last inequality simply uses $\cosh\geqslant1$. This is the required lower bound on the Wiener volume of generic balls. 

\subsection{Identification of a subset with large Wiener volume}

We fix some $\varepsilon>0$ and some $\delta>0$. For reasons that will become clear later, the radius we choose is 
\begin{equation}\label{def:R}R=(\sqrt{2\lambda_d}+\delta)/\varepsilon,
\end{equation} and we finally define an inflated set by $S_0=B(0,\varepsilon)+H_R$. This set will represent the ``high-probability region'' of the Wiener space that will be effectively quantized. 
The celebrated Gaussian isoperimetry theorem gives
\begin{equation}\label{eq:isoperimetry}\mathbb{W}(S_0)\geqslant\Phi\big(\Phi^{-1}(b_d(\varepsilon))+R\big), \end{equation}
where $\Phi(t) = \int_{-\infty}^t \frac{1}{\sqrt{2\pi}}e^{-x^2/2}dx$ is the cumulative distribution function of $\mathscr{N}(0,1)$. By \eqref{eq:bound_bd} and classical asymptotics\footnote{We use that  $\Phi^{-1}(u)\sim -\sqrt{-2\log u}$ as $u\to 0$, and $1-\Phi(u) \leqslant e^{-O(u^2)}$ as $u\to \infty$. } we get that $\Phi^{-1}(b_d(\varepsilon))\sim -\sqrt{2\lambda_d}/\varepsilon$, so that \begin{equation}\label{eq:vol_s0}1 -\mathbb{W}(S_0)\leqslant e^{-O(\varepsilon^{-2})}.\end{equation}

\subsection{Quantization of the large subset}

We start by drawing $m$ auxiliary Brownian paths, say $F_1, \dotsc, F_m$, with
$$m=\left\lceil\frac{\varepsilon^{-2}}{\beta}\right\rceil,\qquad
\beta=e^{-R^2/2}b_d(\varepsilon).$$
These paths are purely auxiliary: they will merely be used to show the existence of an efficient quantization of the Wiener space, but in a non-constructive way. We retain only the paths with a norm smaller than $R+\varepsilon$. The other ones are just thrown away. 

We have seen earlier that if $g$ is a fixed function in $H_R$, then each of these random paths will land in $B(g,\varepsilon)$ with probability at least $\beta$. For any fixed $f$ in $S_0$ there is a $g$ in $H_R$ at distance at most $\varepsilon$: consequently, there is a probability at least $\beta$ that the draw will land in $B(f,2\varepsilon)$, and thus the probability that none of the $m$ paths lands in $B(f,2\varepsilon)$ is at most $e^{-m\beta}\leqslant e^{-\varepsilon^{-2}}$. 

Let us define the total Wiener mass of $S_0$ which is \emph{not} covered by $B(F_i, 2\varepsilon)$: 
$$A=A(F_1, \dotsc, F_m) = \mathbb{W}\left(S_0 \setminus \bigcup_{i=1}^m B(F_i, 2\varepsilon)\right).$$
This is a random variable depending on $F=(F_1, \dotsc, F_m)$. But Fubini's theorem says that 
$$\mathbb{E}[A] = \mathbb{E}_W[\mathbf{1}_{W \in S_0}\mathbb{P}_{F}(\text{none of the } F_i \text{ lands in } B(W, 2\varepsilon))].$$
We have proved that for any fixed $f\in S_0$ (here, we take $f$ to be $W\in S_0$), 
$$\mathbb{P}_{F}(\text{none of the } F_i \text{ lands in } B(W, 2\varepsilon)) \leqslant e^{-\varepsilon^{-2}}.$$
Hence we see that $\mathbb{E}[A] \leqslant e^{-\varepsilon^{-2}}$: so there must exist\footnote{Indeed, this argument shows that there are many such families, but we will not use this fact.} at least one deterministic family with at most $m$ paths, say $f_1, \dotsc, f_m$ for simplicity, such that $\mathbb{W}\left(S_0 \setminus \cup_{i=1}^m B(f_i, 2\varepsilon)\right) \leqslant e^{-\varepsilon^{-2}}$.

\subsection{Definition of the quantization function}
We fix such a family, add the function $f_0=0$, and denote the whole new family by $\mathscr{F}=\{f_0, f_1,\ldots,f_m\}$. All the elements of the family have a norm smaller than $R+\varepsilon$. Let \begin{equation}\label{def:S}S=\bigcup_{j=1}^m B(f_j,2\varepsilon).\end{equation}Let us now define the ``quantization function'' $\mathfrak{q}:\mathscr{C}\to \mathscr{F}$: we map every $f\in\mathscr{C}$ to the closest function among $f_1,\dotsc, f_m$, except if $f$ is not in $S$, in which case we simply set $\mathfrak{q}(f)=f_0=0$. Our choices of $m$ and $R$ entail 
$$\log(m+1)=\frac{(\sqrt{2\lambda_d}+\delta)^2/2+\lambda_d+o(1)}{\varepsilon^2}.$$
For any fixed $\gamma>2\lambda_d$ we may choose $\delta$ small enough so that, for all small $\varepsilon$, the number of elements in the family is $m\leqslant e^{\gamma/\varepsilon^2}$, and the set $S$ defined above has a very large Wiener measure: we have $\mathbb{W}(S^c) \leqslant \mathbb{W}(S_0^c)  + \mathbb{W}(S_0 \setminus \cup_j B(f_j, 2\varepsilon))$, hence thanks to \eqref{eq:vol_s0} we get
$$1 - \mathbb{W}(S) \leqslant 2e^{-O(\varepsilon^{-2})}.$$
Finally, for any fixed function $f \in S$, we have $\Vert f-\mathfrak q(f)\Vert \leqslant 2\varepsilon$. 
This is the only construction needed in both proofs.

\section{Proof of the main theorem in the Brownian case}\label{sec:proof_brownian}

\subsection{Proof of the expected lower bound} \label{sec:brownian_lower}
The proof of the lower bound essentially follows the one in \cite{dereich2005randomly}, bounding the quantity $L_n$ defined below. 
Any coupling between $W \sim \mathbb{W}$ and $Y\sim \mu_n$ has to move the mass to a distance at least equal to the distance between $W$ and the closest $W_i$, so clearly
$$L_n:=\mathbb E_W\left[\min_{i\leqslant n} \Vert W-W_i\Vert^p\right] \leqslant\W_p^p(\mu_n,\mathbb W). $$
Here $\mathbb E_W$ only averages over the new path $W$. Let $E_r=\bigcup_{i=1}^nB(W_i,r)$. Then, 
$$L_n \geqslant \mathbb{E}    [\mathbf{1}_{W \notin E_r}\min_{i\leqslant n} \Vert W-W_i\Vert^p] \geqslant r^p (1 - \mathbb{E}[\mathbb{W}(E_r)]).  $$
The union bound gives $\mathbb{E}[\mathbb{W}(E_r)]\leqslant n b_d(r/\sqrt2)$. Let us fix some $\delta>0$ and take 
\begin{equation}\label{def:r}r^2=\frac{2\lambda_d}{(1+\delta)\log n}.\end{equation} By \eqref{eq:bound_bd}, the expected covered mass is thus $O(n^{-\delta})$, hence
\begin{align*}\mathbb{E}\W_p^p(\mu_n,\mathbb{W})
    &\geqslant r^p(1 + O(n^{-\delta}))\\ 
    &= \left(\frac{2\lambda_d}{(1+\delta)\log n}\right)^{p/2}(1+O(n^{-\delta})).
\end{align*}
Sending $\delta$ to zero proves that $\liminf \sqrt{\log n}(\mathbb{E}\W_p^p(\mu_n,\mathbb{W}))^{1/p}\geqslant \sqrt{2\lambda_d} $. We now must pass to an almost sure result. 

Along $n=2^k$, the discussion above shows that the expected masses $\mathbb{E}[\mathbb{W}(E_r)]=O(n^{-\delta})$ are summable. A simple application of the Borel--Cantelli lemma shows that the mass $\mathbb{W}(E_r)$ tends to zero almost surely, and thus almost surely $L_n\geqslant r_n^p (1+o(1))$ along $n=2^k$. But $L_n$ is decreasing in $n$, so for any $2^k < n \leqslant 2^{k+1}$,
$$\W_p^p(\mu_n,\mathbb{W}) \geqslant L_{2^{k+1}} \geqslant (r_{2^{k+1}})^p(1+o(1)), \quad \text{a.s}.$$
Since $\log 2^{k+1}\sim \log n$ we get that $r_{2^{k+1}} \sim r_n$, as requested. This argument also proves that almost surely, $\sqrt{\log n}L_n$ is asymptotically greater than $\sqrt{2\lambda_d}$.

\subsection{Proof of the upper bound}

Let $\mathbb{W}_\varepsilon$ be the distribution of the quantized Brownian motion $\mathfrak q(W)$, and let
$$\nu_n=\frac1n\sum_{i=1}^n\delta_{\mathfrak q(W_i)}.$$
Both measures are supported on the finite family $\mathscr{F}$, which has $m$ elements, $m$ being smaller than $e^{\gamma\varepsilon^{-2}}$. We therefore fix $0<\alpha<1$ and $\gamma>2\lambda_d$, and we choose $\varepsilon$ so that $m \leqslant n^\alpha$, namely:
$$\varepsilon^2=\frac{\gamma}{\alpha\log n}.$$

\subsubsection{The quantization error}

Write $e(f)=\Vert f-\mathfrak q(f)\Vert$. As we saw earlier, if $W$ lands in $S$, then $e(W)\leqslant 2\varepsilon$. We set
$$D=\mathbb{E}[e(W)^p].$$
This is what we call the \emph{quantization error}. 
By the properties of the quantization function and the Cauchy--Schwarz inequality,
\begin{align*}D&= \mathbb{E}[e(W)^p \mathbf{1}_{W \in S}] + \mathbb{E}[e(W)^p \mathbf{1}_{W \notin S}]\\
    &\leqslant(2\varepsilon)^p+\sqrt{\mathbb{E}[\Vert W\Vert^{2p}]}\sqrt{1-\mathbb{W}(S)}\\
    &\leqslant(2\varepsilon)^p (1 + o(1)).
\end{align*}

\subsubsection{Triangle inequality}
The triangle inequality shows that
\begin{equation}\label{eq:triangle}\W_p(\mu_n,\mathbb{W})\leqslant \W_p(\mu_n, \nu_n) + \W_p(\nu_n, \mathbb{W}_\varepsilon) + \W_p(\mathbb{W}_\varepsilon, \mathbb{W}). \end{equation}

In the definition of $\W_p(\mu_n, \nu_n)$ we can test the simplest possible coupling, which pairs $W_i$ with $\mathfrak{q}(W_i)$. Its cost is $U_n := n^{-1}\sum_{i=1}^n e(W_i)^p$, which is also an upper bound for the Wasserstein distance. Similarly, we can couple $W$ with $\mathfrak{q}(W)$, which gives
$\mathbb{E}\W_p^p(\mathbb{W}_\varepsilon,\mathbb{W})\leqslant D$. Thus, from \eqref{eq:triangle} we get
\begin{equation}\label{eq:three_couplings}
\W_p(\mu_n,\mathbb{W})\leqslant U_n^{1/p}
+\W_p(\nu_n,\mathbb{W}_\varepsilon)+D^{1/p}.
\end{equation}

\subsubsection{The empirical quantized error}

The Voronoï cells of the chosen quantization are defined as $C_j=\{f:\mathfrak q(f)=f_j\}$. The Wiener mass of $C_j$ is precisely the mass assigned by $\mathbb{W}_\varepsilon$ to $\{f_j\}$, and we denote it by $q_j$. The empirical mass assigned to $\{f_j\}$ by $\nu_n$ is $\bar q_j=\frac1n\sum_{i=1}^n\1_{W_i\in C_j}$. 

It is well known that if $\mu$ and $\nu$ are two probability measures supported on a finite set $A$, then $\W_p^p(\mu,\nu)\leqslant \mathrm{diam}(A)^p \times (1/2) \times \sum_{a\in A}|\mu\{a\} - \nu\{a\}|$. Since here the finite family $\mathscr{F}$ is composed of functions whose norm is always smaller than $R+\varepsilon$, we see that $\mathrm{diam}(\mathscr{F})\leqslant 2(R+\varepsilon)$, and thus 
\begin{equation}\label{eq:Vn}\W_p^p(\nu_n,\mathbb{W}_\varepsilon)\leqslant[2(R+\varepsilon)]^p V_n / 2\end{equation} where
$V_n=\sum_{j=0}^m|\bar q_j-q_j|$. 
But the random variable $n\bar q_j$ is nothing but a binomial with parameters $n$ and $q_j$. The Cauchy--Schwarz inequality applied twice gives
$$\mathbb{E}V_n\leqslant\frac1{\sqrt n}\sum_{j=0}^m\sqrt{q_j(1-q_j)}
\leqslant O\left(\sqrt{\frac{m}{n}}\right), $$
and thus 
\begin{equation}\label{eq:middleterm}
    \mathbb{E}[\W_p^p(\nu_n, \mathbb{W}_\varepsilon)] = O\left((R+\varepsilon)^p\sqrt{\frac{m}{n}}\right).
\end{equation}

\subsubsection{Putting it all together}
Taking the $L^p$ norm in \eqref{eq:three_couplings} and then using $\mathbb{E}[U_n]= D$, we obtain
$$(\mathbb{E}\W_p^p(\mu_n,\mathbb{W}))^{1/p}
\leqslant2D^{1/p} + O((R+\varepsilon)(m/n)^{ 1/(2p)})
\leqslant 4\varepsilon(1+o(1)). $$
In particular, $\limsup \sqrt{\log n}(\mathbb{E}\W_p^p(\mu_n,\mathbb{W}))^{1/p} \leqslant 4\sqrt{\gamma/\alpha}$. Letting $\alpha\uparrow1$ and $\gamma\downarrow2\lambda_d$ finishes the proof for the expectation. 

For the almost sure bound, apply Lemma \ref{lem:averages} to the centered variables $e(W_i)^p$ and the centered cell-indicator functions. All their moments are uniformly bounded, since $e(W_i)\leqslant\Vert W_i\Vert + 2\varepsilon$. For every $\eta> 0$, one can use basic concentration inequalities to show that almost surely,
$$|U_n-D| =O(n^\eta /\sqrt{n})\quad\text{ and }\quad V_n=O(n^\eta \sqrt{m/n}).$$
Choose $\eta<(1-\alpha)/2$. So in \eqref{eq:three_couplings} we can replace $U_n$ by $D$ up to a polynomially small error, while \eqref{eq:Vn} gives a polynomially small bound for the middle term. So directly from \eqref{eq:triangle} we get $\W_p(\mu_n,\mathbb{W})\leqslant4\varepsilon(1+o(1))$ almost surely, and we conclude as above by taking the limsup and then $\alpha\uparrow1$, $\gamma \downarrow2\lambda_d$.

\section{Proof of the main theorem in the Girsanov case}\label{sec:proof_girsanov}

We now explain how to adapt the preceding proof to the Girsanov case. The main quantization argument is the same, but now we will need to control the weights $Z_i$ and their normalization. We recall that
$$S_n = \frac{1}{n}\sum_{i=1}^n Z_i.$$

\subsection{Preliminaries}

Let $K=\sup_{t,x}|h(t,x)|$. For $q\geqslant1$ or $q<0$, compare $Z^q$ with $\tilde{Z}$, the exponential martingale \eqref{eq:girsanov_weight} but with integrand $qh$: their ratio is smaller than $e^{q(q-1)K^2/2}$, so
$$\mathbb{E}[Z_1^q] = \mathbb{E}[\tilde{Z}(Z_1^q/\tilde{Z})] \leqslant \mathbb{E}[\tilde{Z}]e^{q(q-1)K^2/2}=e^{q(q-1)K^2/2}.$$
Consequently, 
\begin{equation}\label{eq:weight_moments}
\mathbb{E}[Z_1^q]\leqslant e^{q(q-1)K^2/2}.
\end{equation}
Also, Jensen's inequality yields that for any $q>0$ we have $\mathbb{E}[S_n^{-q}]\leqslant\mathbb{E}[Z_1^{-q}]$. 

\subsection{Proof of the lower bound}

We keep the notation $E_r=\bigcup_{i=1}^n B(W_i,r)$, and define $Z$ from $W$ by \eqref{eq:girsanov_weight}. Just as for the unweighted case, the mass assigned by $\mathbb{Q}$ to the complement of $E_r$ has to be transported by a distance of at least $r$. But then, the Girsanov theorem \eqref{eq:girsanov_expectation}, the Cauchy--Schwarz inequality and the bound on the moments of $Z$ yield
\begin{align*}
    \mathbb{E}[\mathbb{Q}(E_r)] &= \mathbb{E}[\mathbb{E}_W[Z \mathbf{1}_{W \in E_r}]]\\
    &\leqslant \sqrt{\mathbb{E}[Z^2]\mathbb{P}(W \in E_r)} \\
    &\leqslant e^{K^2/2}\sqrt{nb_d(r/\sqrt{2})}.
\end{align*}
With $r^2=2\lambda_d/((1+\delta)\log n)$ as before, this tends to zero and is summable along $n=2^k$. The proof in Section~\ref{sec:brownian_lower} applies unchanged, provided we use $\mathbb{Q}$ in the definition of $L_n$. In particular, the leading constant does not depend on $K$.

\subsection{Proof of the upper bound}
We still use the quantizer $\mathfrak q$, with $\varepsilon^2=\gamma/(\alpha\log n)$ where $\gamma>2\lambda_d$ and $0<\alpha<1$. In particular, $m\leqslant n^\alpha$. Let $\mathbb{Q}_\varepsilon$ be the law of $\mathfrak q(X)$, and let
$$\hat{\nu}_n=\frac{1}{nS_n}\sum_{i=1}^n Z_i\delta_{\mathfrak q(W_i)}.$$
We reproduce the strategy that led to \eqref{eq:three_couplings}, that is, we first use the triangle inequality
$$\W_p(\widehat\mu_n, \mathbb{Q}) \leqslant
\underbrace{\W_p(\widehat\mu_n, \hat{\nu}_n)}_{\text{(I)}}
+\underbrace{\W_p(\hat{\nu}_n, \mathbb{Q}_\varepsilon)}_{\text{(II)}}
+\underbrace{\W_p(\mathbb{Q}_\varepsilon, \mathbb{Q})}_{\text{(III)}}$$
and then we bound each term separately, following the Brownian proof with adaptations.

\subsubsection{Bound on $\text{(I)}$}

We set
$$U_n=\frac{1}{n}\sum_{i=1}^nZ_i e(W_i)^p.$$
Clearly $\mathfrak{q}$ gives a coupling between $\widehat\mu_n$ and $\hat{\nu}_n$, hence $\text{(I)} \leqslant (U_n/S_n)^{1/p}$. 
It should be emphasized that both measures are normalized to have mass 1. We split the cost between $S$ and its complement, then use exchangeability:
$$\mathbb{E}[U_n/S_n]\leqslant(2\varepsilon)^p+
\mathbb{E}\left[\frac{Z_1\Vert W_1\Vert^p}{S_n}\1_{W_1\notin S}\right].$$
To bound the last term, we use the multi-term H\"older inequality with exponents $8,8,4,2$:
$$\mathbb{E}\left[\frac{Z_1\Vert W_1\Vert^p}{S_n}\1_{W_1\notin S}\right]
\leqslant\mathbb{E}[Z_1^8]^{1/8}\mathbb{E}[\Vert W_1\Vert^{8p}]^{1/8}\,\mathbb{E}\big[S_n^{-4}\big]^{1/4}\,\mathbb{W}(S^c)^{1/2}.$$
Thanks to the bounds on the moments of $Z$ and $S_n$ at the beginning of this section, the first term is bounded, and the third one too. Hence the first three terms are bounded by a constant which depends neither on $n$ nor on $\varepsilon$. The last term is clearly $o(\varepsilon^p)$, and we can conclude as in the Brownian case.

\subsubsection{Bound on $\text{(II)}$}

Let us denote the masses assigned by $\mathbb{Q}_\varepsilon$ and $\hat\nu_n$ to the singleton $f_j$ by (respectively) $q_j=\mathbb{E}[Z\1_{W\in C_j}]$ and $\bar{q}_j =A_j/S_n$, where $ A_j=\frac{1}{n}\sum_{i=1}^n Z_i\1_{W_i\in C_j}$. Here again, 
$$\W_p^p(\mathbb{Q}_\varepsilon, \hat{\nu}_n)\leqslant O((R+\varepsilon)^p V_n)$$
where $V_n = \sum_j |q_j - \bar{q}_j|$. 
We have
\begin{equation}\label{eq:weighted_imbalance}
V_n  \leqslant |S_n-1|+\sum_j|A_j-q_j|.
\end{equation}
Indeed, $|A_j/S_n-q_j|\leqslant|A_j/S_n-A_j|+|A_j-q_j|$ by the triangle inequality, and since $A_j\geqslant0$ we have $\sum_j|A_j/S_n-A_j|=|1-S_n|\sum_jA_j/S_n=|1-S_n|$.

Now, let us write $v=\mathbb{E}Z^2\leqslant e^{K^2}$. By independence, $\mathbb{E}\sum_j(A_j-q_j)^2=(v-\sum_jq_j^2)/n = O(1/n)$ and $\mathbb{E}(S_n-1)^2=(v-1)/n = O(1/n)$. Thus the Cauchy--Schwarz inequality yields $\mathbb{E}V_n = O\left(\sqrt{\frac{m}{n}}\right)$. 
This recovers the Brownian bound. 

\subsubsection{Bound on $\text{(III)}$}
Since $\mathfrak{q}$ induces a coupling between $\mathbb{Q}$ and $\mathbb{Q}_\varepsilon$, we have $\W_p^p(\mathbb{Q}_\varepsilon, \mathbb{Q}) \leqslant D:=\mathbb{E}[e(X)^p]$. Thanks to Girsanov's theorem,  $D =\mathbb{E}[Ze(W)^p]$. 
As before, the Cauchy--Schwarz inequality (applied twice) and the bounds on the moments of $Z$ give
$$D\leqslant(2\varepsilon)^p+\sqrt{\mathbb{E}[Z^2\Vert W\Vert^{2p}](1-\mathbb{W}(S))}\leqslant(2\varepsilon)^p(1+o(1)).$$
This is exactly the same bound as for the Brownian case.

\subsection{Conclusion}
Gathering the bounds obtained for the three terms $\text{(I), (II), (III)}$, we get $$\big(\mathbb{E}\W_p^p(\widehat\mu_n,\mathbb{Q})\big)^{1/p}
\leqslant D^{1/p}+\big(\mathbb{E}[U_n/S_n]\big)^{1/p}
+2(R+\varepsilon)(\mathbb{E}V_n)^{1/p}
\leqslant4\varepsilon(1+o(1)).$$

This proves the bound in expectation. The ultimate missing step is the almost-sure bound, which is done in the same way. Indeed, apply Lemma \ref{lem:averages} to the centered summands of $U_n$, the centered vectors $(Z_i\1_{W_i\in C_j})_j$, and $Z_i-1$. Their moments are uniformly bounded. For every $\eta>0$, almost surely,
$$|U_n-D|+|S_n-1|=O(n^{-1/2+\eta}),\qquad
\left(\sum_j(A_j-q_j)^2\right)^{1/2}=O(n^{-1/2+\eta}).$$
Choose $\eta<(1-\alpha)/2$. By \eqref{eq:weighted_imbalance}, $V_n=O(n^{-(1-\alpha)/2+\eta})$. Since $S_n\to1$, the three couplings yield
$$\W_p(\widehat\mu_n,\mathbb{Q})\leqslant D^{1/p}+ (U_n /S_n)^{1/p}+O((R+\varepsilon)V_n^{1/p})
\leqslant4\varepsilon(1+o(1)).$$
Letting $\alpha\uparrow1$ and $\gamma\downarrow2\lambda_d$ as before finishes the proof.

\bibliographystyle{abbrvurl}
\bibliography{OT}

@article{chatterjee_diaconis,
  author = {Sourav Chatterjee and Persi Diaconis},
  title = {{The sample size required in importance sampling}},
  volume = {28},
  journal = {The Annals of Applied Probability},
  number = {2},
  publisher = {Institute of Mathematical Statistics},
  pages = {1099 -- 1135},
  year = {2018},
  doi = {10.1214/17-AAP1326},
  URL = {https://doi.org/10.1214/17-AAP1326}
  }

@article{goldman2024optimal,
  title={Optimal transport methods for combinatorial optimization over two random point sets},
  author={Goldman, Michael and Trevisan, Dario},
  journal={Probability Theory and Related Fields},
  volume={188},
  number={3},
  pages={1315--1384},
  year={2024},
  publisher={Springer}
}

@article{agapiou2017importance,
  title={Importance sampling: Intrinsic dimension and computational cost},
  author={Agapiou, Sergios and Papaspiliopoulos, Omiros and Sanz-Alonso, Daniel and Stuart, Andrew M},
  journal={Statistical Science},
  pages={405--431},
  year={2017},
  publisher={JSTOR}
}

@article{kahn1953methods,
  title={Methods of reducing sample size in Monte Carlo computations},
  author={Kahn, Herman and Marshall, Andy W},
  journal={Journal of the Operations Research Society of America},
  volume={1},
  number={5},
  pages={263--278},
  year={1953},
  publisher={INFORMS}
}

@article{elvira2021advances,
  title={Advances in importance sampling},
  author={Elvira, V{\'\i}ctor and Martino, Luca},
  journal={arXiv preprint arXiv:2102.05407},
  year={2021}
}

@InCollection{BaBo13,
    Author = {F. {Barthe} and C. {Bordenave}},
    Title = {{Combinatorial optimization over two random point sets.}},
    BookTitle = {{S\'eminaire de probabilit\'es XLV}},
    ISBN = {978-3-319-00320-7/pbk; 978-3-319-00321-4/ebook},
    Pages = {483--535},
    Year = {2013},
    Publisher = {Cham: Springer},
    DOI = {10.1007/978-3-319-00321-4_19},
    MSC2010 = {90C27}
}

@Article{AKT84,
    Author = {M. {Ajtai} and J. {Koml\'os} and G. {Tusn\'ady}},
    Title = {{On optimal matchings.}},
    FJournal = {{Combinatorica}},
    Journal = {{Combinatorica}},
    ISSN = {0209-9683; 1439-6912/e},
    Volume = {4},
    Pages = {259--264},
    Year = {1984},
    Publisher = {Springer, Berlin/Heidelberg; J\'anos Bolyai Mathematical Society, Budapest},
    DOI = {10.1007/BF02579135},
    MSC2010 = {60D05 60G40},
    Zbl = {0562.60012}
}

@article {FoGu15,
    AUTHOR = {Fournier, N. and Guillin, A.},
     TITLE = {On the rate of convergence in {W}asserstein distance of the
              empirical measure},
   JOURNAL = {Probab. Theory Related Fields},
  FJOURNAL = {Probability Theory and Related Fields},
    VOLUME = {162},
      YEAR = {2015},
    NUMBER = {3-4},
     PAGES = {707--738},
      ISSN = {0178-8051},
   MRCLASS = {60F25 (60E15 60F10)},
MRREVIEWER = {Jos\'e Trashorras},
       DOI = {10.1007/s00440-014-0583-7},
       URL = {https://doi.org/10.1007/s00440-014-0583-7},
}

@book{graf2000foundations,
  title={Foundations of quantization for probability distributions},
  author={Graf, Siegfried and Luschgy, Harald},
  year={2000},
  publisher={Springer Science \& Business Media}
}

@Article{Ta94,
    Author = {M. {Talagrand}},
    Title = {{The transportation cost from the uniform measure to the empirical measure in dimension $\geq 3$.}},
    FJournal = {{The Annals of Probability}},
    Journal = {{Ann. Probab.}},
    ISSN = {0091-1798; 2168-894X/e},
    Volume = {22},
    Number = {2},
    Pages = {919--959},
    Year = {1994},
    Publisher = {Institute of Mathematical Statistics (IMS), Beachwood, OH/Bethesda, MD},
    DOI = {10.1214/aop/1176988735},
    MSC2010 = {60D05 60C05},
    Zbl = {0809.60015}
}

@article{boissard2014mean,
  title={On the mean speed of convergence of empirical and occupation measures in {W}asserstein distance},
  author={Boissard, Emmanuel and Le Gouic, Thibaut},
  journal={Ann. Inst. Henri Poincar\'e Probab. Stat.},
  volume={50},
  number={2},
  pages={539--563},
  year={2014},
  doi={10.1214/12-AIHP517}
}

@misc{coste2023importance,
  author={Coste, Simon},
  howpublished={Blog post},
  url={https://scoste.fr/posts/importance/}
}

@misc{ghrepo,
  author={Coste, Simon},
  howpublished={GitLab repository},
  url={https://gitlab.com/SimonCoste/wasserstein}
}

@article{coste2026wasserstein,
  title={The {W}asserstein cost of importance sampling},
  author={Coste, Simon and Goldman, Michael},
  journal={arXiv preprint arXiv:2605.30055},
  year={2026}
}

@article{dereich2005randomly,
  author={Dereich, Steffen and Lifshits, Mikhail A.},
  title={Probabilities of randomly centered small balls and quantization in {B}anach spaces},
  journal={Ann. Probab.},
  volume={33},
  number={4},
  pages={1397--1421},
  year={2005},
  doi={10.1214/009117905000000161}
}

@article{dereich2003link,
  author={Dereich, Steffen and Fehringer, Franz and Matoussi, Anis and Scheutzow, Michael},
  title={On the link between small ball probabilities and the quantization problem for {G}aussian measures on {B}anach spaces},
  journal={J. Theoret. Probab.},
  volume={16},
  number={1},
  pages={249--265},
  year={2003},
  doi={10.1023/A:1022242924198}
}

@article{dereich2006high,
  author={Dereich, Steffen and Scheutzow, Michael},
  title={High-resolution quantization and entropy coding for fractional {B}rownian motion},
  journal={Electron. J. Probab.},
  volume={11},
  number={28},
  pages={700--722},
  year={2006}
}

@article{dereich2008coding,
  author={Dereich, Steffen},
  title={The coding complexity of diffusion processes under supremum norm distortion},
  journal={Stochastic Process. Appl.},
  volume={118},
  number={6},
  pages={917--937},
  year={2008},
  doi={10.1016/j.spa.2007.07.003}
}

@article{newton1994variance,
  author={Newton, Nigel J.},
  title={Variance reduction for simulated diffusions},
  journal={SIAM J. Appl. Math.},
  volume={54},
  number={6},
  pages={1780--1805},
  year={1994},
  doi={10.1137/S0036139992236220}
}

@article{glasserman1999asymptotically,
  author={Glasserman, Paul and Heidelberger, Philip and Shahabuddin, Perwez},
  title={Asymptotically optimal importance sampling and stratification for pricing path-dependent options},
  journal={Math. Finance},
  volume={9},
  number={2},
  pages={117--152},
  year={1999},
  doi={10.1111/1467-9965.00065}
}

@book{glasserman2004monte,
  author={Glasserman, Paul},
  title={Monte {C}arlo methods in financial engineering},
  series={Applications of Mathematics},
  volume={53},
  publisher={Springer},
  address={New York},
  year={2004}
}

@article{he2026sharp,
  title={Sharp functional quantization and empirical {W}asserstein rates for {I}t\^o processes},
  author={He, Xihao and Yuan, Fengyi},
  journal={arXiv preprint arXiv:2608.07879},
  year={2026}
}

@article{weed2019sharp,
  title={Sharp asymptotic and finite-sample rates of convergence of empirical measures in Wasserstein distance},
  author={Weed, Jonathan and Bach, Francis},
  journal={Bernoulli},
  volume={25},
  number={4A},
  pages={2620--2648},
  year={2019},
  publisher={JSTOR}
}

@article{anderson1955integral,
  author={Anderson, T. W.},
  title={The integral of a symmetric unimodal function over a symmetric convex set and some probability inequalities},
  journal={Proc. Amer. Math. Soc.},
  volume={6},
  number={2},
  pages={170--176},
  year={1955}
}

@article{ciesielski1962first,
  author={Ciesielski, Z. and Taylor, S. J.},
  title={First passage times and sojourn times for {Brownian} motion in space and the exact {Hausdorff} measure of the sample path},
  journal={Trans. Amer. Math. Soc.},
  volume={103},
  number={3},
  pages={434--450},
  year={1962},
  doi={10.1090/S0002-9947-1962-0143257-8}
}

\appendix

\section{Small-ball probabilities}\label{appendix:small_balls}

\subsection{The one-dimensional formula}

In dimension one, let $u(t,x)$ be the probability that $x+W$ stays in $]-r,r[$ up to time $t$. It solves the heat equation with boundary conditions:
$$\partial_tu=\frac{1}{2}\partial_{xx}u,\qquad
u(t,\pm r)=0,\qquad u(0,x)=1.$$
Expanding the constant function in the Dirichlet eigenfunctions and evaluating at $t=1,x=0$ gives the exact formula for $b_1$:
\begin{equation}\label{eq:small_ball_series}
    b_1(r)=\frac4\pi\sum_{j\geqslant0}\frac{(-1)^j}{2j+1}
    \exp\left(-\frac{(2j+1)^2\pi^2}{8r^2}\right).
\end{equation}
This is an alternating series whose terms are decreasing in absolute value, so the first term bounds the series from above. For the lower bound, the first eigenfunction $\cos(\pi x/(2r))$ is between zero and one on $]-r,r[$ and equals 1 at the origin. Its killed heat evolution gives
$$e^{-\pi^2/(8r^2)}\leqslant b_1(r)\leqslant\frac4\pi e^{-\pi^2/(8r^2)}.$$
Both prefactors are optimal among bounds of this form, by letting $r\to\infty$ and $r\downarrow0$, respectively. Independence of the coordinates gives the elementary bounds
$$b_1(r/\sqrt d)^d\leqslant b_d(r)\leqslant b_1(r)^d,\qquad
 e^{-\pi^2d^2/(8r^2)}\leqslant b_d(r)\leqslant(4/\pi)^d e^{-\pi^2d/(8r^2)}.$$

\subsection{The constant in dimension $d$}

Brownian scaling identifies $b_d(r)$ with survival in the unit ball up to time $r^{-2}$. The Dirichlet heat expansion has a positive first coefficient, so
$-\log b_d(r)=\lambda_d/r^2+O(1)$. Solving the radial eigenvalue equation gives
$$\lambda_d=\frac12j_{d/2-1,1}^{\,2},\qquad
\lambda_1=\frac{\pi^2}{8},\quad \lambda_2\simeq2.8916,\quad
\lambda_3=\frac{\pi^2}{2},$$
where $j_{\nu,1}$ is the first positive zero of the Bessel function $J_\nu$; see \cite{ciesielski1962first}. This proves \eqref{eq:bound_bd}. 

\subsection{Anderson's inequality}

For a centered Gaussian vector, the probability of a translation of a symmetric convex set is largest when the translation is zero \cite{anderson1955integral}. Apply this to the vector $(W(t_k))_k$ on a time grid and to the product of Euclidean balls of radius $r$. For any fixed path $f$,
$$\mathbb{P}\left(\max_k|W_1(t_k)-f(t_k)| \leqslant r\right) \leqslant\mathbb{P}\left(\max_k|W_1(t_k)|\leqslant r\right).$$
Letting the grids increase to a dense subset of $[0,1]$ and using continuity gives \eqref{eq:anderson}.

\section{An elementary averaging estimate}\label{appendix:averaging}

\begin{lemma}\label{lem:averages}
For each $n$, let $Y_{n,1}, \ldots,Y_{n,n}$ be independent, identically distributed, centered Euclidean vectors. Their dimension may depend on $n$. If
$\sup_n\mathbb{E}\Vert Y_{n,1}\Vert_2^{2k}<\infty$ for every integer $k\geqslant1$, then for every $\eta>0$, almost surely one has
$$\left\Vert\frac1n\sum_{i=1}^nY_{n,i}\right\Vert_2=O(n^{-1/2+\eta}).$$
\end{lemma}

\begin{proof}
Expand the $k$th power of $\sum_{i,j}\langle Y_{n,i},Y_{n,j}\rangle$. A term with an index appearing only once has expectation zero. The remaining terms involve at most $k$ distinct indices, so there are $O_k(n^k)$ of them; H\"older's inequality bounds each term by the $2k$-th moment. Thus 
$$\mathbb{E} \left\Vert\frac1n\sum_{i=1}^n Y_{n,i}\right\Vert_2^{2k}\leqslant C_k n^{-k}\mathbb{E}\Vert Y_{n,1}\Vert_2^{ 2k}. $$
By Markov's inequality the probability that the average is larger than $n^{-1/2+\eta}$ is bounded by $C_k n^{-2k\eta}$. Set $2k\eta>1$ and apply the Borel--Cantelli lemma to conclude. 
\end{proof}
   
\noindent \hrulefill

\end{document}